\documentclass[11pt, a4paper]{amsart}
\usepackage{graphicx} 

\usepackage[letterpaper,top=2cm,bottom=2cm,left=3cm,right=3cm,marginparwidth=1.75cm]{geometry}
\usepackage{xcolor}
\definecolor{myblue}{rgb}{.8, .8, 1}
\usepackage{amsmath}
\usepackage{amsfonts,amssymb}
\usepackage[pagewise]{lineno}\linenumbers
\usepackage{empheq}

\newcommand{\IS}{\mathbb{S}}
\newcommand{\R}{\mathbb{R}}
\newcommand{\tr}{\mathrm{Tr}}
\newcommand{\Tr}{\mathrm{Tr}}
\newcommand{\td}{\text{d}}
\newcommand{\ito}{It\^o~}
\newcommand{\Ito}{\ito}

\newcommand{\T}{\mathbb{T}}

\newlength\mytemplen
\newsavebox\mytempbox

\makeatletter
\newcommand\mybluebox{%
    \@ifnextchar[
       {\@mybluebox}%
       {\@mybluebox[0pt]}}

\def\@mybluebox[#1]{%
    \@ifnextchar[
       {\@@mybluebox[#1]}%
       {\@@mybluebox[#1][0pt]}}

\def\@@mybluebox[#1][#2]#3{
    \sbox\mytempbox{#3}%
    \mytemplen\ht\mytempbox
    \advance\mytemplen #1\relax
    \ht\mytempbox\mytemplen
    \mytemplen\dp\mytempbox
    \advance\mytemplen #2\relax
    \dp\mytempbox\mytemplen
    \colorbox{myblue}{\hspace{1em}\usebox{\mytempbox}\hspace{1em}}}

\usepackage{ bbold }
\usepackage{hyperref}
\usepackage{color}
\newcommand{\sagy}[1]{{\color{blue}{[Sagy: #1]}}}
\newcommand{\erik}[1]{{\color{orange}{[Erik: #1]}}}
\newcommand{\andrea}[1]{{\color{red}{[Andrea: #1]}}}

\title{Diffusive behavior of transport noise on $\IS^2$}
\author{Sagy Ephrati}
\author{Erik Jansson}
\author{Andrea Papini}
\address[Sagy Ephrati, Erik Jansson, Andrea Papini]{Department of Mathematical Sciences, Chalmers University of Technology and University of Gothenburg, 412 96 Gothenburg, Sweden}
\email{\{sagy, erikjans, andreapa\}@chalmers.se}
\date{\today}
\usepackage{cleveref}
\usepackage{comment}
\begin{document}

\maketitle
\begin{abstract}
    We investigate theoretically and numerically transport noise-induced diffusion in flows on the sphere. Previous analysis on the torus demonstrated that suitably chosen transport noise in the Euler equations leads to diffusive behavior resembling the Navier--Stokes equations. Here, we analyze dynamics on the sphere with noise-induced differential elliptic operator dissipation and characterize their energy and enstrophy decay properties. Through structure-preserving numerical simulations with the Zeitlin discretization, we demonstrate that appropriately scaled transport noise induces energy dissipation while preserving enstrophy and coadjoint orbits. The presented analysis lays a groundwork for further theoretical investigation of transport noise and supports the calibration of transport noise models as a parametrization for unresolved processes in geophysical fluid simulations.
\end{abstract}
\section{Introduction}
The broad range of length and time scales in fluid-dynamical systems presents a major challenge in practical applications where accurate flow predictions are essential. 
To make computations tractable, physical models are often simplified, and numerical simulations are truncated at finite resolutions. 
This creates a need for additional closure models that account for unresolved dynamics and mitigate the effects of numerical approximations. 
An emerging approach within geophysical fluid dynamics is the use of \textit{transport noise} \cite{holm2015variational, memin2014fluid} on coarse computational grids as a surrogate for computationally expensive high-resolution simulations. 
This framework offers a structure-preserving parametrization of unresolvable scales and model uncertainty, enabling physically consistent stochastic flow predictions at reduced computational costs.
Numerical studies have demonstrated the effectiveness of transport noise for efficient uncertainty quantification in advection-dominated flows \cite{cotter2019numerically, ephrati2023data, resseguier2017geophysical} with promising applications in data assimilation \cite{cotter2020particle, dufee2022stochastic}.

\begin{figure}[h]
\centering
 \includegraphics[width = 0.9\linewidth]{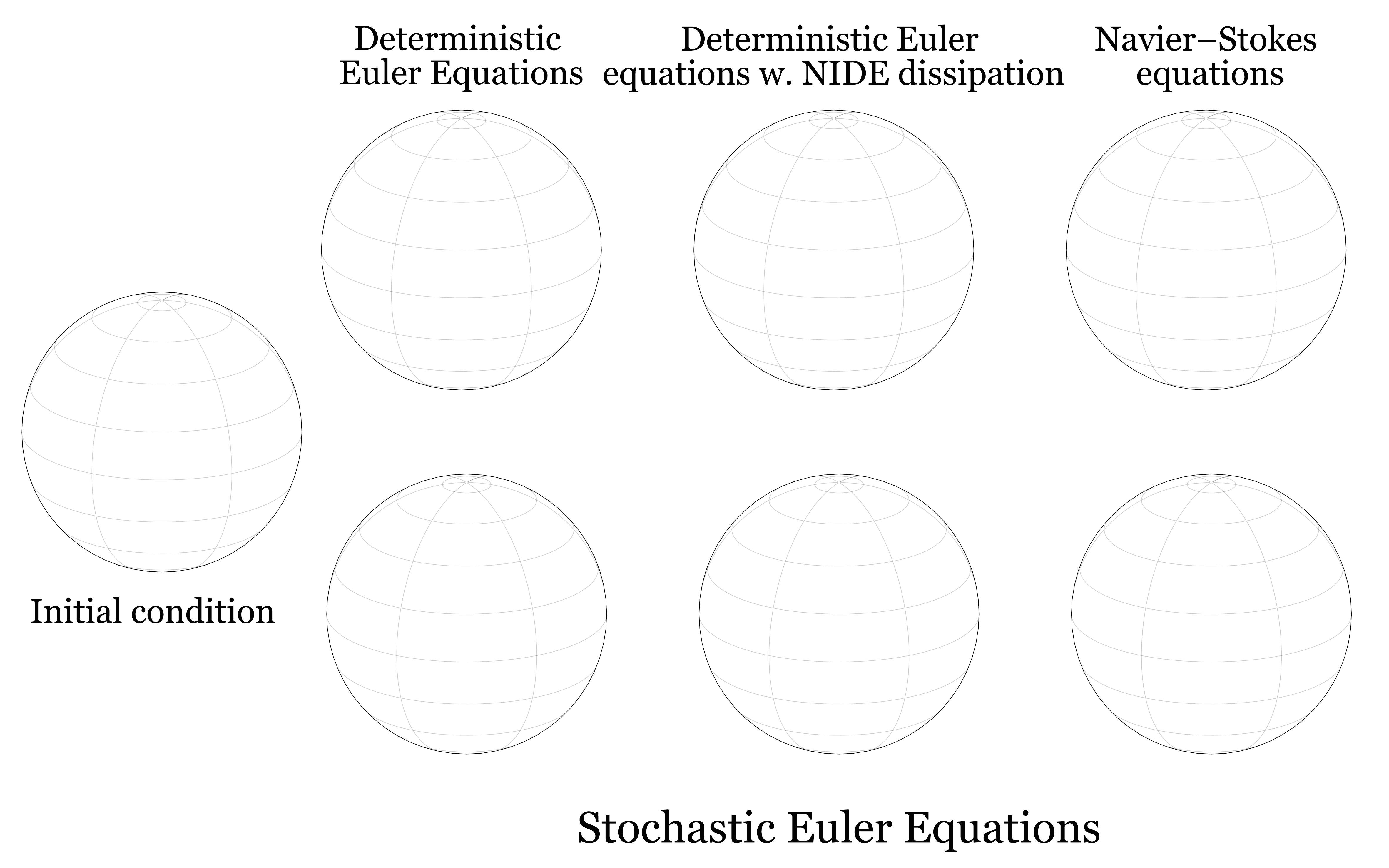}
 \caption{Snapshots of fluid flow on the sphere, achieved from the same initial condition but with different governing equations. Deterministic equations are used in the top row, stochastic equations with various noise scalings are used in the bottom row.}
 \label{fig:eyecandy}
\end{figure}

The increased use of transport noise in computational studies has spurred analytical developments to aid the calibration and understanding of stochastic models \cite{lang2023well, goodair20223d}.
Our work is motivated by the seminal paper \cite{flandoli2021scaling}, where it was proven that a suitable choice of transport noise in the Euler equations on the torus leads to diffusive behavior resembling the Navier--Stokes equations. In particular, it is shown in \cite{flandoli2021scaling} that the correct choice of transport noise properties results in a dissipation of energy. 
The goal of this paper is to investigate noise-induced dissipation by structure-preserving simulation. 
The enabling factor of the current study is the development of a scalable geometric integrator for the Euler equations on the sphere subject to transport noise \cite{ephrati2024exponential}. 
An example of flow on the sphere is given in \Cref{fig:eyecandy}.

The theoretical results by \cite{flandoli2021scaling} unite the favorable geometric properties of transport noise with dissipative properties of traditionally used reduced-complexity closure models. 
However, the resemblance between the stochastic Euler equations and the Navier--Stokes equations arises in a limiting regime. 
The proof is delicate and places precise requirements on the noise to achieve the necessary convergence properties. 
On the torus, this is achieved by constructing the noise as a linear combination of Fourier modes. 
Extending these results to other manifolds poses a challenge, since these manifolds generally lack explicitly known Fourier bases, thus requiring redefinition of the noise. 
A promising alternative is the sphere, where one can use spherical harmonic functions to define distinct noise structures. Nonetheless, the exact way to converge to Navier--Stokes via transport noise on the sphere is currently unknown, and a rigorous proof is expected to require substantial theoretical machinery. This is to be investigated carefully in future work.

This paper has two main objectives. 
First, it lays the groundwork for a detailed analysis of transport noise in spherical geometries.
Second, we perform structure-preserving simulations that illustrate diffusive behavior of transport noise to support the theoretical findings of \cite{flandoli2021scaling}. 

\section{Euler equations: Deterministic, stochastic, and dissipative}
We introduce the two-dimensional Euler equations on the sphere and highlight some of their geometry. 
Subsequently, the stochastic Euler equations and dissipative Euler equations are introduced, which are central to the present work.

\subsection{Deterministic Euler equations}
The deterministic Euler equations on $\IS^2$ are given by 
\begin{align}
    \label{eq:euler}
    \partial_t \omega + \{ \psi, \omega \} = 0, \quad \Delta_{\IS^2} \psi = \omega 
\end{align}
where $\omega$ is the vorticity, $\psi$ is the stream function, and $\Delta_{\IS^2}$ is the Laplace--Beltrami operator.
The Poisson bracket is defined as
\begin{align}
    \{ f, g \} = \nabla^\perp f \cdot \nabla g,
\end{align}
where $\nabla^\perp$ is the skew-gradient obtained by rotating the gradient by 90 degrees. 

The Euler equations form a Lie--Poisson system on the Lie--Poisson algebra $\{C^\infty_0(\IS^2),\{\cdot,\cdot\}\}$ and has several conserved quantities. 
First, the energy is conserved, given by
\begin{align} \label{eq:energy}
    E(\omega) = -\frac{1}{2} \int_{\IS^2} \psi\omega \, d\mu,
\end{align}
where $d\mu$ is the volume form on $\IS^2$.
Second, there is an infinite number of conserved quantities called \emph{Casimirs}. 
For any smooth function $f: \IS^2 \to \R$, the corresponding Casimir is given by
\begin{align}
    C_f(\omega) = \int_{\IS^2} f(\omega) \, d\mu.
\end{align}
A well-known example of a Casimir is the \emph{enstrophy}, given by 
\begin{align*}
    S(\omega) = \int_{\IS^2} \omega^2 \, d\mu.
\end{align*}

Due to the Lie--Poisson structure of the Euler equations, solutions of the Euler equations evolve on subspaces known as \emph{coadjoint orbits}. 
This distinct geometric structure of the Euler equations is conjectured to have important consequences for the long-term behavior. 
Moreover, the geometric properties of the Euler equations are not shared by the Navier--Stokes equations, 
\begin{align*}
    \partial_t \omega + \{ \psi, \omega \} + \nu \Delta \omega = 0,
    \label{eq:NS}
\end{align*}
where $\nu$ is the viscosity.
In particular, the Navier--Stokes equation does not have any conserved quantities, nor does the solution evolve on coadjoint orbits.

\subsection{Euler equations with transport noise}
As outlined above, the Euler equations have a distinct geometric structure. 
If one wishes to make the equations stochastic while preserving this structure, one must introduce stochasticity in the form of \emph{transport noise} to \Cref{eq:euler}. 
This yields
\begin{equation}\label{EqSph}
    \td\omega + \{\psi, \omega \}\td t+\sum_{i = 1}^M \{\psi_i,\omega\} \circ \td B_t^i=0,
\end{equation}
where $\psi_i:\IS^2\rightarrow\R$ are scalar-valued smooth functions 
and $B_t^1,\ldots,B_t^N$ are $N$ independent scalar Brownian motions.
The multiplication here is understood in the Stratonovich sense. The Stratonovich integral is used in \eqref{EqSph}, since the ordinary chain rule holds for Stratonovich processes \cite{kloeden1992stochastic}. 
In turn, this facilitates the extension of manifold-valued curves in differential geometry to manifold-valued processes \cite{emery2006two}. 
The solutions to the stochastic Euler equations \eqref{EqSph} still evolve on the coadjoint orbits and have the Casimirs as conserved quantities. 
However, the energy \eqref{eq:energy} is not necessarily conserved. 

The Stratonovich formulation of \Cref{EqSph} provides a geometric way to include stochasticity in the Euler equations.
However, in this paper, we are interested in dissipative properties of transport noise, for example obtained as limiting behavior as the number of noise terms goes to infinity. 
For this reason, it is more useful to consider a formulation that consists of a deterministic term and one zero-mean stochastic term, i.e., to move to the It\^o formulation. 
We employ the It\^o correction \cite{holm2015variational}
to obtain the equation 
\begin{equation}\label{EqSphIto}
    \td\omega + \{\psi, \omega \}\td t - \frac{1}{2}\sum_{i = 1}^M\{\psi_i,\{\psi_i,\omega\}\}\td t+\sum_{i = 1}^M \{\psi_i,\omega\} \td B_t^i=0. 
\end{equation}
\Cref{EqSphIto} now contains two additional terms compared to the deterministic Euler equations \eqref{eq:euler}. 
First, there is a zero-mean martingale term 
\begin{align}
    \label{eq:noise_operator}
    N(\omega) = \sum_{i = 1}^M \{\psi_i,\omega\} \td B_t^i
\end{align}
and second, a dissipative term 
\begin{align}
    \label{eq:NIDE_operator}
     \Lambda_M \omega = \frac{1}{2}\sum_{i = 1}^M\{\psi_i,\{\psi_i,\omega\}\}.
\end{align}
We refer to the operator $\Lambda_M $ as the \emph{noise-induced differential elliptic (NIDE) operator}, which depends on the choice of $\psi_i$.

Adding any one of the terms \eqref{eq:noise_operator} or \eqref{eq:NIDE_operator}, but not the other, results in a drift away from the coadjoint orbits. 
The idea of \cite{flandoli2021scaling} is to prove that the $L^2$-norm of $N_M(\omega)$ goes to $0$ as $M \to \infty$, and that $\Lambda_M$ in turn converges, in an appropriate norm, to a non-zero operator.  
The resulting system does not, thus, in the limit, evolve on coadjoint orbits. In fact, in the case of $\T^2$ when the sequence $(\psi_\ell)_{\ell = 0}^\infty$ is chosen appropriately as in \cite{flandoli2021scaling}, the NIDE operator $\Lambda_M$ converges to the toroidal Laplacian in a suitable sense. Thereby, in the limit, solutions to the Navier--Stokes equations are recovered in the sense of $H^{-1}$ convergence.

We remark, however, that for any choice of finite $M$, a  choice of non-zero functions $(\psi_\ell)_{\ell = 0}^M$ results in a \emph{structure-preserving} system, in the sense that \Cref{EqSph} still evolves on the coadjoint orbits. This is not the case for, e.g., the Navier--Stokes equations. In the scaling limit considered by \cite{flandoli2021scaling}, the stochastic system converges, possibly in a suitable norm, to a solution that no longer preserves the geometric structure present with finitely many noise terms. While each system with a finite number of terms evolves on a coadjoint orbit, the limit appears to lie outside this solution space, indicating a breakdown or deformation of the underlying geometry rather than a structure-preserving convergence.
Further, not all choices of $(\psi_\ell)_{\ell = 0}^M$ result in energy dissipation, and energy conservation can be achieved with certain selections of constant or solution-dependent scalar fields. 
A particular example is the case where $\psi_1 = \psi = \Delta^{-1}\omega$ and $\psi_i=0$ for $i=2,\ldots, M$. 
The resulting system is thus given by 
\begin{align}
    \label{eq:secret_av}
    \td \omega + \{\psi,\omega\}\left(\mathrm d t + \circ \mathrm B_t^1 \right) = 0.
\end{align}
Assuming the right solution regularity, we have the following energy evolution:
\begin{align*}
    \td E &=-\left[\int_{\IS^2} \{\psi,\omega\}\psi  d\mu\right] (\td t+ \circ \td B_t^1) = 0, 
\end{align*}
since $\int_{\IS^2} \{\psi,\omega\}\psi = -\int_{\IS^2} \{\psi,\psi\omega\} = 0$. 
Hence, we see that the energy is preserved using this choice of NIDE-operator.

If we instead select the scalar fields analogous to the  $\T^2$--case in \cite{flandoli2021scaling}, 
we obtain the sequence 
\begin{align}
    \label{eq:francochoice}
    \begin{split}
    \psi^M_{\ell,m}&=\alpha^M_{\ell,m} Y_{\ell,m}\\
    \alpha^M_{\ell,m}&=\frac{\sqrt{2\nu_\mathrm{SALT}}}{\|\{c_{k,n}\}_{k=1,n =-k,\ldots,k}^{k=\infty}\|_{\ell^2}}c_{\ell,m},\quad c^M_{\ell,m}=\frac{1}{(\ell+1)^a}\mathbb{1}_{\{0\leq\ell\leq M\}},
    \end{split}
\end{align} 
where $l \geq 0$ and $m = -l,\ldots, l$.
Here, the parameter $a>0$ prescribes the possible convergence behavior of the two operators \eqref{eq:noise_operator}, \eqref{eq:NIDE_operator} according to the correct Sobolev norm. 
We expect to retrieve dissipative behavior from such a selection since $\{\psi_{\ell,m},\omega\}$ is not zero everywhere and the selection of $\psi_{l,m}$ does not depend on the solution $\omega$.
We will show results for the convergence of \eqref{EqSphIto} to a deterministic dissipative system, using heuristic computations. Through the limiting NIDE operator, although it is \emph{not} the Laplace--Beltrami operator on the sphere, we obtain an elliptic dissipative operator. 
Several numerical experiments are carried out to support this.

\subsection{Euler equations with NIDE operator dissipation}
To better understand the dissipative behavior of the stochastic Euler equations, we study the deterministic equations with the dissipative term appearing in the \Ito formulation \eqref{EqSphIto}. We refer to these as the \emph{Euler equations with NIDE operator dissipation}, which are given by 
\begin{align}
    \label{eq:nide_euler}
    \td\omega + \{\psi, \omega \}\td t - \Lambda_M \omega = 0.
\end{align}
The system \eqref{eq:nide_euler} contains several important fluid dynamical equations. 
Indeed,  by setting $\psi_1 = Y_{1,0}$, $\psi_2 = Y_{1,-1}$ and $\psi_3 = Y_{1,1}$, and all other $\psi_\ell$ to $0$, 
we have that by \cite[Section 3.1]{modin2024two}
\begin{align*}
    \Lambda_M \omega = \sum_{i = 0}^3 \{Y_{1,i},\{Y_{1,i},\omega\}\} \propto \Delta \omega. 
\end{align*}
Thus, the Laplacian is a NIDE operator. 
Another NIDE operator appears in the anticipated vorticity method (AVM), and is obtained with $\psi_1 = \sqrt{2}\psi = \Delta^{-1} \omega$ and all other functions set to zero,  
\begin{align}
    \label{eq:nide_av}
    \Lambda_M \omega = \{\psi,\{\psi,\omega\}\}.
\end{align}
In the deterministic setting, this operator decays the enstrophy while leaving the energy unaffected, resulting in a type of \emph{selective decay} \cite{Sadourny1985, GayBalmaz2013}.
The NIDE operator \eqref{eq:nide_av} is in fact the It\^o--Stratonovich corrector appearing for the system \eqref{eq:secret_av}, for which we observed previously that the energy remained constant. 
In contrast to the Euler equations with NIDE operator dissipation given by \eqref{eq:nide_av}, the system \eqref{eq:secret_av} remains on coadjoint orbits, meaning that the It\^o noise term compensates for the enstrophy dissipation. 


Concerning the properties of the NIDE operator for a selection of sufficiently smooth functions $\phi_i$, $\Lambda_M$ can dissipate energy as well as enstrophy. 
To this end, consider the example system 
\begin{align*}
    \dot f = \Lambda_M f,
\end{align*}
where $f$ is in the domain of $\Lambda_M$.
To investigate the enstrophy dissipation properties, we compute
\begin{align*}
    \frac{\mathrm{d}}{\mathrm{d}t} S(f) &=\int_{\IS^2} \frac{1}{2}\sum_{\ell = 1}^M \{\psi_\ell,\{\psi_\ell,f\}\}f\td \mu=\frac{1}{2}\sum_{\ell = 1}^M \int_{\IS^2} \{\psi_\ell,\{\psi_\ell,f\}\}f\td \mu\\
    &=\frac{1}{2}\sum_{\ell = 1}^M \int_{\IS^2} \left(f\nabla^\perp\psi_\ell\right)\cdot\nabla\left(\{\psi_\ell,f\}\right)\td \mu\\
    &=-\frac{1}{2}\sum_{\ell = 1}^M \int_{\IS^2} \left(\nabla f\cdot\nabla^\perp\psi_\ell\right)\{\psi_\ell,f\}\td \mu-\frac{1}{2}\sum_{\ell = 1}^M \int_{\IS^2} f\text{div}\left(\nabla^\perp\psi_\ell\right)\{\psi_\ell,f\}\td \mu\\
    &=-\frac{1}{2}\sum_{\ell = 1}^M \int_{\IS^2} \{\psi_\ell,f\}^2\td \mu\leq0.
\end{align*}
Here, we have used the divergence theorem on the sphere and the fact that $\nabla^\perp\psi_\ell$ is a divergence--free vector field for every choice of $\psi_\ell$. 

Likewise, for the energy, we find
\begin{align*}
    \frac{\mathrm{d}}{\mathrm{d}t} E(f) 
    &=-\frac{1}{2}\sum_{\ell = 1}^M \int_{\IS^2} \{\psi_\ell,\Delta^{-1} f\} \{\psi_\ell, f\}\td \mu = -\frac{1}{2}\sum_{\ell = 1}^M \int_{\IS^2} \{\Delta^{-1} \psi_\ell, f\} \{\psi_\ell, f\}\td \mu. 
\end{align*}
This computation shows that the dissipative behavior of $\Lambda_M$ depends on how $\psi_\ell$ interacts with $\Delta^{-1}$. 
For instance, if $\psi_\ell = \Delta^{-1} f$ as in the AVM operator \eqref{eq:nide_av}, then $\frac{\mathrm{d}}{\mathrm{d}t} E(f) 
     = 0$. 
    
However, if we insert the sequence in \Cref{eq:francochoice} and use that $\Delta^{-1} \alpha_{\ell,m} Y_{\ell,m} = \frac{\alpha_\ell}{\ell(\ell+1)} Y_{\ell,m} $, we obtain 
\begin{align}
    \label{eq:energy_diss}
    \begin{split}
    \frac{\mathrm{d}}{\mathrm{d}t} E(f) 
    &=-\frac{1}{2}\sum_{\ell = 1}^M\sum_{m = -\ell}^\ell  \int_{\IS^2} \{\psi_\ell,\Delta^{-1} f\} \{\psi_{\ell,m}, f\}\td \mu \\
    &= -\frac{\nu_\mathrm{SALT}}{{\|\{c_{k,n}\}_{k=1,n =-k,\ldots,k}^{k=\infty}\|_{\ell^2}}^2} \sum_{\ell = 1}^M \sum_{m = -\ell}^\ell \frac{1}{\ell(\ell+1)^{2a+1}} \int_{\IS^2} \{Y_{\ell,m}, f\}^2\td \mu \leq 0.
    \end{split}
\end{align}
This means that energy dissipates when using the sequence in \Cref{eq:francochoice}. 
Similarly, the enstrophy dissipates as 
\begin{align}
    \label{eq:enstrophy_diss}
    \frac{\mathrm{d}}{\mathrm{d}t} S(f)  = 
    & -\frac{\nu_\mathrm{SALT}}{{\|\{c_{k,n}\}_{k=1,n =-k,\ldots,k}^{k=M}\|_{\ell^2}}^2}\sum_{\ell = 1}^M \sum_{m=-\ell}^\ell\frac{1}{(\ell+1)^{2a}} \int_{\IS^2} \{Y_{\ell,m}, f\}^2\td \mu \leq 0.
\end{align}
Clearly, enstrophy is dissipated faster than energy. 

These results illustrate the effect of the noise scaling parameter $a$ on the rate of dissipation. Specifically, for each $M>0$, the total dissipation depends on $a\in \mathbb{R}_+$ in such a way that increasing $a$ places more emphasis on the low-degree modes in the NIDE operator.
Namely, as $a$ grows, the quantity $\|c_{\cdot}\|_{\ell^2}^{-2}(\ell+1)^{-2a-1}$ reduces the magnitudes corresponding to high-degree modes in the NIDE-operator. 
Thus, the noise scaling has significant consequences. If chosen too large, the Boussinesq hypothesis which states that small-scale turbulent features act on large-scale features in a dissipating manner, becomes irrelevant and the scaling limit of the stochastic equation is not guaranteed. This line of reasoning inspired \cite{flandoli2021scaling}. On the other hand, if $a$ is too small, then the NIDE operator induces slow dissipation and a greater deviation from the Navier--Stokes is expected.




\subsection{Necessary requirement for dissipative noise limit on the sphere}

We now analyze briefly the choice of noise parameters \eqref{eq:francochoice}, the possible scaling limit and the behavior of the noise term.

Let us consider solutions $\omega^M$ and $\omega$ of \eqref{EqSphIto} and \eqref{eq:nide_euler}, respectively.
Postponing the technical details to future work, we conjecture that $\omega^M\rightarrow\omega$ in a suitably weak sense and that in particular, given $\phi\in C^\infty(\mathbb S^2)$,
\begin{align}
    \mathbb{E}\left[\left<\omega^M(t)-\omega(t),\phi\right>^2\right]&\lesssim 
    \mathbb{E}\left[\left<\omega^M(0)-\omega(0),\phi\right>^2\right]+\mathbb{E}\left[\left<\int_0^t\{\psi^M,\omega^M(s)\}\td s-\int_0^t\{\psi,\omega(s)\}\td s,\phi\right>^2\right]\nonumber\\
&\quad+\mathbb{E}\left[\left<\sum_{\ell=0}^{M}\sum_{m=-\ell}^\ell(\alpha^M_{\ell,m})^2\nabla^\perp Y_{\ell,m}\cdot\nabla(\nabla^\perp Y_{\ell,m}\cdot \nabla\phi),\omega\right>^2\right]+ \mathcal{M}_t^M.\label{eq:weakOp-noise}
\end{align}
Here, the term $\mathcal{M}_t^M$ represents the stochastic part of \eqref{EqSphIto}, 
$$
\mathcal{M}^M_t:=\mathbb{E}\left[\left|\sum_{\ell=0}^{M}\sum_{m=-\ell}^\ell\int_0^t\left<(\alpha^M_{\ell,m})\nabla^\perp Y_{\ell,m}\cdot \nabla\phi),\omega\right>\ dB_t^{\ell,m}\right|^2\right],
$$
and is only affected by $\psi, \alpha^M$ and $\omega^M$.
It is evident that convergence of the solution of the stochastic equation to the solution of the Euler equations with NIDE operator dissipation at a minimum requires that $\mathcal{M}^M\rightarrow0_{M\rightarrow\infty}$. 
Let us then consider the following:
\begin{align*}
\mathcal M^M_t&\leq\mathbb{E}\left[\sup_{t\in[0,T]}\left|\sum_{\ell=0}^M\sum_{m=-\ell}^\ell \int_0^t \alpha_{\ell,m}^M\left<\nabla^\perp Y_{\ell,m}^M\cdot \nabla \phi,\omega\right>\td B_t^{\ell,m}\right|^2\right]\\
&\leq \mathbb{E}\left[\sum_{\ell=0}^M\sum_{m=-\ell}^\ell \int_0^T (\alpha_{\ell,m}^M)^2\left<\nabla^\perp Y_{\ell,m}^M\cdot \nabla \phi,\omega\right>^2\td t\right]\\
&\leq \|\alpha_\cdot^M\|_\infty^2\sum_{\ell=0}^M\sum_{m=-\ell}^\ell\mathbb{E}\left[ \int_0^T \left<\nabla^\perp Y_{\ell,m}^M\cdot \nabla \phi,\omega\right>^2\td t\right]\leq\|\alpha_\cdot^M\|_\infty^2\mathbb{E}\left[ \int_0^T \|\omega\nabla\phi\|_{L^2}^2\td t\right]\\
&\leq \|\alpha_\cdot^M\|_\infty^2 \|\nabla\phi\|_\infty\mathbb{E}\left[ \int_0^T \|\omega\|_{L^2}^2\td t\right].
\end{align*}
Hence, if we select the parameter $a\in \mathbb R_+$ in \eqref{eq:francochoice} such that $\|\alpha_\cdot^M\|_\infty\rightarrow0$, the noise vanishes.

For this reason, looking at \eqref{eq:francochoice}, we note that
$$
\|\alpha^M_\cdot\|^2_{L^2}=2\nu_\mathrm{SALT}\sum_{\ell=0}^M\sum_{m=-\ell}^\ell c_{\ell,m}^2\|\{c_{\ell,m}\}\|_{\ell^2}^{-2}=2\nu_\mathrm{SALT}.
$$
This shows that the first term in \eqref{eq:weakOp-noise} is bounded since we can bound each term with $\mathbb E\|\omega\|_{L^2}^2$ and a constant depending on the gradient of the spherical harmonics and the test function $\phi$, leaving the $\ell^2$-norm of our selected sequence.
At the same time, we obtain the $\ell^\infty$-norm of this sequence as
\begin{align*}
    \|\alpha^M_\cdot\|_\infty=\sup_{\ell,m}\|\{c_{\ell,m}\}\|_{\ell^2}^{-1}\frac{\mathbb{1}_{\{0\leq\ell\leq M\}}}{(\ell+1)^a}=\|\{c_{\ell,m}\}\|_{\ell^2}^{-1},
\end{align*}
since the supremum is taken over all possible combinations of $\ell,m$ with $0\leq \ell\leq M$, and we have that $\sup_{\{0\leq\ell\leq M, m\}} \frac{1}{(l+1)^a}=1$ for every $a\in \R_+$.
Using estimates on the generalized harmonic series, we find that asymptotically
\begin{align*}
    \|\{c_{\ell,m}\}\|_{\ell^2}&\sim\begin{cases}
        (M+1)^{a-1},& 0\leq a<1,\\
        \left(\log(M+1)\right)^{-\frac{1}{2}}, & a=1,\\
        1/\sqrt{2},& a>1.\\
    \end{cases}
\end{align*}
Hence, when $M\rightarrow\infty$ $\|,\{\alpha^M_{\ell,m}\}\|_{\ell^\infty}$ converges to zero only when $0\leq a\leq1$. It is only possible in this case, at least in the distributional sense, for solutions of the stochastic equation \eqref{EqSphIto} to converge to those of the deterministic NIDE--equation \eqref{eq:nide_euler}.

Thus, when $a>1$, we observe that the NIDE operator appears to regularize the solutions by curbing high-frequency components similar to a Laplacian as seen in \eqref{eq:energy_diss}; however, convergence fails due to the constraint $0\leq a \leq 1$ when attempting to pass to the limit. Consequently, taking the dissipative limit is expected to also affect the solution at intermediate frequencies.

\section{Computational method}
The dissipative properties of transport noise are supported by numerical simulations, presented in the subsequent section. 
In this section, we introduce the geometric integration methods used to obtain reliable numerical results.
Using a structure-preserving discretization is crucial in the present study, as it guarantees that the numerical solutions remain on coadjoint orbits and that adverse discretization effects are minimal. 
We adopt the Euler--Zeitlin equations as a spatial discretization for the Euler equations on the sphere, which are integrated in time with a stochastic Lie--Poisson integrator.
We highlight these methods below and provide a discussion on the interpretation of noise-induced dissipation in discretized systems.

\subsection{Discretized equations}
We adopt the self-consistent finite-dimensional truncation method of Zeitlin \cite{zeitlin2004self}, for incompressible flow on the sphere as a spatial discretization for the governing equations.
The method provides a finite-dimensional approximation through geometric quantization \cite{hoppe1989diffeomorphism, bordemann1991gl, bordemann1994toeplitz}, via a projection $\Pi_N:\mathcal{C}^\infty(S^2)\to \mathfrak{u}(N)$ of smooth functions on the sphere to complex $N\times N$ skew-Hermitian matrices. Here, $N$ can be regarded as the numerical resolution.
This process replaces the Poisson bracket $\{\cdot, \cdot\}$ by the matrix commutator $\frac{1}{\hbar}\left[\cdot,\cdot\right]$ where $\hbar=2/\sqrt{N^2-1}$. The resulting discretized dynamics read
\begin{equation}
    \dot{W} + \frac{1}{\hbar}[P, W] = 0, \quad \Delta_N P = W
    \label{eq:Zeitlin_Euler}
\end{equation}
where $W$ is the vorticity matrix, and where $P$ is the stream matrix. We will refer to this system as the \textit{Zeitlin--Euler equations}.
The discretized system has a Lie--Poisson structure and possesses Casimir functions of the form 
\begin{equation}
    C_f = \tr(f(W)).
\end{equation} 
The solutions of the Zeitlin--Euler equations \eqref{eq:Zeitlin_Euler} evolve on coadjoint orbits, analogous to those of the Euler equation \eqref{eq:euler}. 
For an initial condition $W_0 = W(0)$, these solutions spaces are given by $\mathcal{O}_{W_0}=\{ g^* W_0 (g^*)^{-1}|g\in SU(N)\}$.

The vorticity matrix and the stream matrix are related via the Hoppe--Yau Laplacian $\Delta_N$ \cite{hoppe1998some}. The latter is a discrete counterpart of the Laplace--Beltrami operator, and has eigen-matrices $T_{lm}$ satisfying $\Delta_N T_{lm}=-l(l+1)T_{lm}$ analogous to the continuous spherical harmonic eigenfunctions $Y_{lm}$ of $\Delta_{\IS^2}$. The analogy between the eigenfunctions of $\Delta_N$ and $\Delta_{\IS^2}$ is used to define the projection operator $\Pi_N$. It acts on a scalar field $\omega$ by truncating the spherical harmonic expansion at degree $N-1$ and replacing the basis functions $Y_{lm}$ by the discrete approximations $T_{lm}$, \begin{equation}
    W = \Pi_N(\omega) = \sum_{l=0}^{N-1}\sum_{m=-l}^{l}\omega_{lm}(iT_{lm}),
\end{equation}
where $\omega_{lm}$ are the spherical harmonic coefficients of $\omega$ and $i$ is the imaginary unit. We refer to \cite{modin2024two} for a detailed description of the Zeitlin--Euler equations. 

The discretization for the Navier--Stokes equations follows readily from the Zeitlin--Euler equations \eqref{eq:Zeitlin_Euler} and the discrete Laplacian. 
Namely, the Hoppe--Yau Laplacian is used to construct a viscous dissipation term $\nu\Delta_NW$, where $\nu$ is the viscosity.
The desired discretization is then obtained by adding this to the right-hand side of \Cref{eq:Zeitlin_Euler}.

Transport noise is included in the Zeitlin--Euler equations \eqref{eq:Zeitlin_Euler} analogously to the continuous case. The added stochastic terms will be expressed as a linear combination of the basis matrices $T_{lm}$, so that we obtain \begin{equation} 
    \td W + \frac{1}{\hbar}[P, W]\td t + \frac{1}{\hbar}\sum_{l=1}^{N-1}\sum_{m=-l}^l[\alpha_{lm}T_{lm}, W] \circ\td B_{t}^{lm}=0.
    \label{eq:Stochastic_Zeitlin_Euler}
\end{equation}
This stochastic system of equations is isospectral and consequently shares the geometric properties of the deterministic system \eqref{eq:Zeitlin_Euler}. That is, the solution space remains unchanged under the addition of transport noise, and therefore the Casimir functions are retained. A detailed geometric description of the deterministic and stochastic system is provided in \cite{ephrati2024exponential}.

The Stratonovich integral is used in \eqref{eq:Stochastic_Zeitlin_Euler}, analogous to the stochastic equations in the continuum \eqref{EqSph}. However, the \ito formulation is well suited for analysis. A general derivation of the \ito form for transport noise is found in \cite{holm2015variational}. Below, we provide the specific form corresponding to \eqref{eq:Stochastic_Zeitlin_Euler}, where we denote the double sum over $l, m$ by $\sum_{lm}$ for readability. We have
\begin{equation}\label{eq:StochasticZeitlinEulerIto}
    \td W + \frac{1}{\hbar}[P, W]\td t + \frac{1}{\hbar}\sum_{lm}[\alpha_{lm}T_{lm}, W]\,\td B_{t}^{lm} = \frac{1}{2}\frac{1}{\hbar^2}\left[\sum_{lm}\alpha_{lm}T_{lm}, \left[\sum_{l'm'}\alpha_{l'm'}T_{l'm'}, W\right]\right]\td t.
\end{equation}

The corresponding Zeitlin--Euler equations with NIDE operator dissipation are given by \begin{equation}\label{eq:ZeitlinEulerNIDE}
    \dot W + \frac{1}{\hbar}[P, W] = \frac{1}{2}\frac{1}{\hbar^2}\left[\sum_{lm}\alpha_{lm}T_{lm}, \left[\sum_{l'm'}\alpha_{l'm'}T_{l'm'}, W\right]\right].
\end{equation}

\subsection{Time integration}
A Casimir-preserving symplectic midpoint method \cite{ephrati2024exponential} is employed for time integration of the isospectral parts of the stochastic and deterministic systems to ensure preservation of the Lie--Poisson structure. We introduce a fixed time step size $h$ and use subscripts $n, n+1$ to indicate the $n^\mathrm{th}$ and ensuing time step. The integration method reads
\begin{equation} \label{eq:stoch_isomp}
    \begin{split}
        \tilde P_{h, n} &= h~\Delta_N^{-1}W + \sqrt{h}~\sum_{lm}\alpha_{lm}T_{lm}(\zeta_{lm})_n \\
        W_n &= \left(I - \frac{1}{2}\tilde P_{h, n}\right) \tilde W \left(I + \frac{1}{2}\tilde P_{h, n}\right) \\
        W_{n+1}&= \left(I + \frac{1}{2}\tilde P_{h, n}\right) \tilde W \left(I-\frac{1}{2}\tilde P_{h, n}\right)
    \end{split}
\end{equation}

Here, $(\zeta_{lm})_n$ are truncated random variables suited for implicit stochastic integration \cite{milstein2004stochastic}. Namely, if $\xi\sim\mathcal{N}(0,1)$, then the truncated variable $\zeta$ is given by \begin{equation}
    \zeta = \begin{cases}
        A & \text{if}\quad \xi > A, \\
        -A & \text{if}\quad\xi < -A, \\
        \xi & \text{if}\quad| \xi |\leq A, \\
    \end{cases}
\end{equation}
where $A = \sqrt{2l|\log h|}$ and $l$ is an arbitrary fixed positive integer, here chosen as $l=2$. These truncated variables are independent at each time step and for each $l$ and $m$. This method has a strong order of convergence $1/2$ and a weak order of convergence of 1.
In the absence of stochastic terms, the integrator \eqref{eq:stoch_isomp} coincides with the second-order isospectral midpoint method for deterministic Lie--Poisson systems \cite{modin2020lie} which is well-suited for simulating ideal spherical hydrodynamics \cite{modin2020casimir}.

\subsection{Discussion on numerical investigation of dissipative limits}

Recent developments in scaling limits of transport noise are predominantly theoretical, see e.g., \cite{flandoli2021scaling, FlaLuo1, carigi2023dissipation, agresti2024anomalous}.  
The appearance of a dissipative term in the stochastic Euler equations as the number of stochastic terms increases presents several challenges for numerical studies, which we highlight below.



First, viscous dissipation does not preserve the geometric structure of the Euler equations. 
Equivalently stated, solutions to the Navier--Stokes equations do not evolve on the coadjoint orbits that characterize the solution spaces of the Euler equations. 
The Zeitlin--Euler equations have a similar geometric structure as the Euler equations, being a finite-dimensional Lie--Poisson system with solutions confined to coadjoint orbits even when incorporating transport noise. This raises the question of how well solutions outside the coadjoint orbits, due to viscous dissipation, can be approximated by solutions constrained to the coadjoint orbits.

Second, and more importantly, any discretization of the stochastic Euler equations is necessarily finite-dimensional, and this limits the number of stochastic terms that can be achieved numerically.
In practice, this number is truncated at some finite (possibly resolution--dependent) value, which implies that the findings by \cite{flandoli2021scaling} are generally not directly transferable to numerical studies using the Zeitlin model. Instead, the exact deterministic dissipative limit can only be approximated. One example can be seen in \cite{flandoli2023average}, where averaging the pathwise solutions can more easily show some behavior of the limiting equation, but only in simplifying settings. While an alternative to this is, for instance, point-vortex dynamics as studied in \cite{FlaPapMor} where diffusion behavior via transport noise is shown, restoring stability to Kelvin--Helmholtz effects.

To elaborate on this point, let $\omega^M$ denote the solution to the stochastic Euler equations \eqref{EqSph} on the 2-torus $\mathbb T^2$ with $M$ noise terms, under the selection of stream functions from \cite{flandoli2021scaling}, and let $\omega_\mathrm{NS}$ denote the solution to the {deterministic 2D Navier--Stokes equations}. 
Flandoli et al. \cite{flandoli2021scaling} establishes, for equivalent in law processes, using Skorokhod embedding, that
\[
\lim_{M \to \infty} \mathbb{E} \| \omega_\mathrm{NS} - \omega^M \|_{C([0,T],H^{-\delta}(\mathbb{T}^2))}^2 = 0,
\]
 in the limit as $M \to \infty$,
That is, $\omega_M$ converges to $\omega_\mathrm{NS}$ over $\mathbb{T}^2$ in the distributional sense for every $\delta>0$.
This result \textit{does not} extend directly to the Zeitlin--Euler truncation $\omega^M$, even for large $M$ or increased numerical resolution.
Thus, analogous limiting behavior for $\omega^M$ cannot be inferred.
In fact, even though the vorticity processes $\omega^M$ converge to a limit
$\omega$ that is explicitly dissipative and should suggest that a partial dissipation should be visible at the level of the stochastic solution $\omega^M$, this is not necessarily the case and such behavior clearly depends on the choices of $\{\psi_i\}$.
Given the solutions $\omega^M$ with almost surely preserved $L^2$-norm for every $M>0$, i.e.,
\[
\mathbb{P} \left( \|\omega^M_t\|_{L^2} = \|\omega^M_0\|_{L^2},\ \forall t \in [0,T] \right) = 1,
\]
the limit $\omega = \lim_{M \to \infty} \omega_M$ corresponds to a \textit{dissipative} solution of the Navier--Stokes equations, where kinetic energy and enstrophy decay over time.
Therefore, convergence in the sense of measure of 
\[
\omega_M \to \omega \quad \text{in} \quad C([0,T]; H^{-})
\]
is \textit{optimal}, as shown in \cite{flandoli2021scaling}. Namely, stronger convergence (e.g., in $C([0,T]; L^2)$) would imply conservation of the $L^2$ norm in the limit, contradicting the known dissipative nature of Navier--Stokes solutions. Several constructions validate this behavior:
\begin{itemize}
    \item In \cite{flandoli2021scaling}, solutions $\omega^M$ conserve vorticity norms but converge in $H^-$ to dissipative deterministic solutions.
    \item If $\omega_0 \in C^\infty(\mathbb{T}^2)$ and $\alpha^M$ are finitely supported for each $M$, then $\omega^M$ remains regular due to bounds on $\|\omega_t^M\|_{L^\infty}$, see \cite{crisan2019solution, flandoli2021scaling}. This leads to
    \[
    \|\omega_t^M\|_{L^2} = \|\omega_0\|_{L^2} \quad \forall t > 0,
    \]
    yet $\omega_M$ still converges weakly to a dissipative limit.
\end{itemize}
The weak convergence $\omega^M \to \omega$ in $C([0,T]; H^-)$ is insufficient to imply rigorously anomalous dissipation, i.e. the enhancement of diffusion at the level of paths for every choice of $\{\psi_i\}$ even when no previous dissipation is present in the system. Dissipation appears in the limit, but not all approximating sequences exhibit dissipative behavior, e.g. the choice \eqref{eq:francochoice} with $a=2$ (presented in the next section), or in the anticipated vorticity method \eqref{eq:secret_av} where energy is preserved every $M$. Consequently, no universal conclusion regarding the dissipative behavior and its amplification can be drawn from the scaling limit alone, hence it is impossible to rigorously show this convergence in the truncated numerical system.


\section{Numerical experiments}
We proceed with ensemble simulations of the stochastic Zeitlin--Euler equations. 
The noise scaling is chosen as \begin{equation}
    \begin{split} \label{eq:noise_scaling_numerical}
        c_{lm} &= \frac{1}{(l+1)^a},\quad l=1,\ldots, M, \\
        \alpha_{lm}&= \frac{\sqrt{2\nu_\mathrm{SALT}}}{\|\{c_{lm}\}\|_{l^2}} c_{lm},
    \end{split}
\end{equation}
where $M$ describes the maximal spherical harmonic degree at which the stochastic forcing is activated and $a$ determines the noise strength as the degree increases. 
All simulations are started from the same initial condition, which is chosen by sampling the first $10$ spherical harmonic coefficients independently from the standard normal distribution.

Both values are user-specified and will be varied in a range of numerical tests. Specifically, we adopt $a=1.0$, analogous to the scaling proposed on the torus \cite{flandoli2021scaling}, and $a=2$, for comparison. 
Note that the $l^2$-norm of the sequence $\{c_{lm}\}$ depends on the values of $M$ and $a$.

The numerical simulations are performed at a resolution $N=256$ with $M=N/16,~N/8, ~N/4, ~N/2$.
The maximum adopted value of $M$ is set to $N/2$ to minimize aliasing effects associated with finite truncation in spectral methods.
The stochastic ensembles are compared to two distinct deterministic systems. 
First, the standard Zeitlin--Navier--Stokes equations are used as a guideline for how viscous dissipation affects the flow. 
Second, the Zeitlin--Euler equations with NIDE operator dissipation \eqref{eq:ZeitlinEulerNIDE} serve as an additional reference for the dissipative behavior in stochastic systems.

The dissipation strength is quantified via the total energy, which for the Zeitlin--Euler equations reads \begin{equation}
    H = -\frac{1}{2}\Tr(PW).
\end{equation}
We additionally measure the enstrophy \begin{equation}
    S = \Tr(WW).
\end{equation} 
The enstrophy is conserved in the standard deterministic Euler equations and the stochastic Euler equations, as a consequence of the solutions evolving on coadjoint orbits.
The enstrophy is generally not conserved for the Navier--Stokes equations and the Euler equations with NIDE operator dissipation.

\subsection{Results}
\begin{figure}[h]
    \centering
    \includegraphics[width=0.49\linewidth]{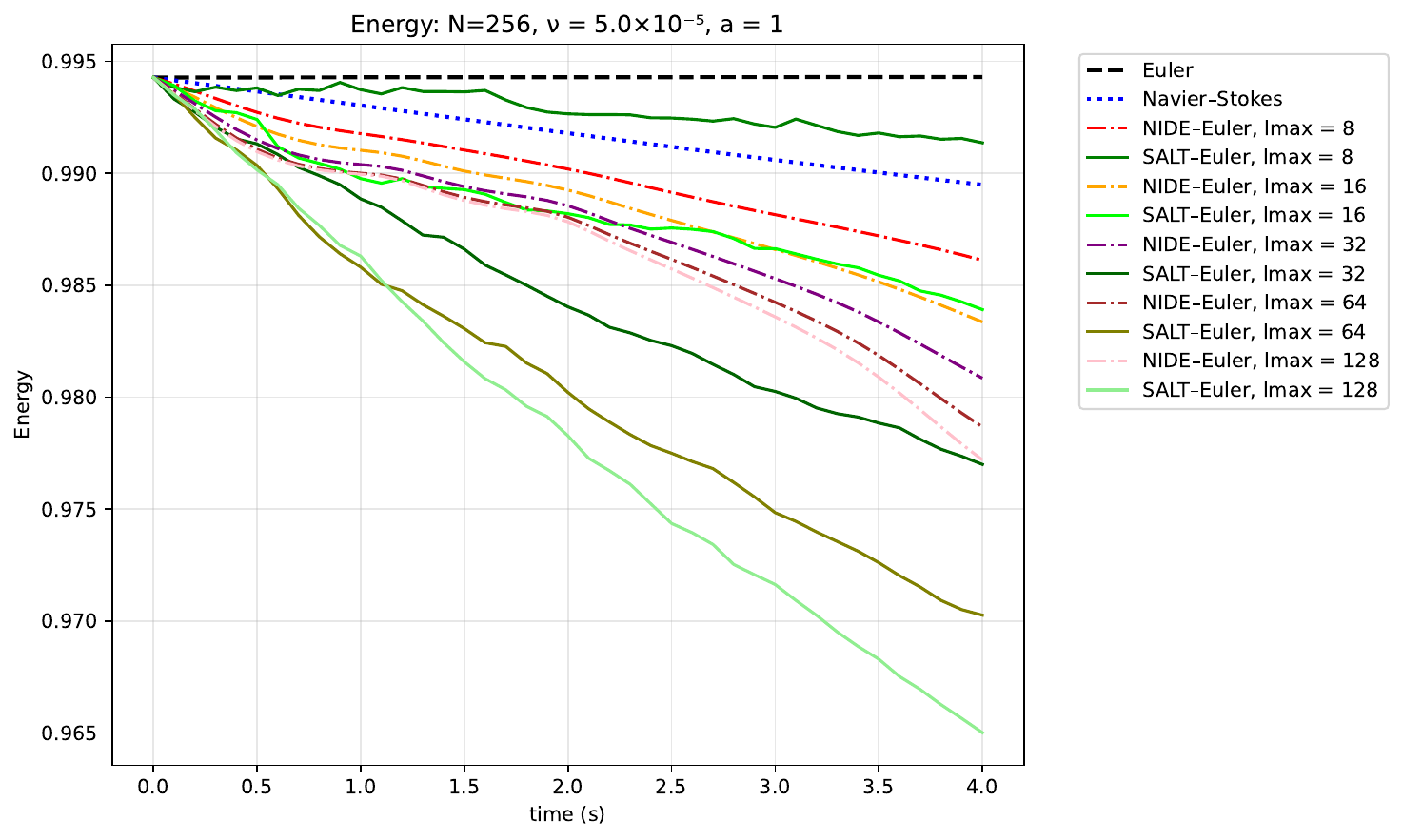}
    \includegraphics[width=0.49\linewidth]{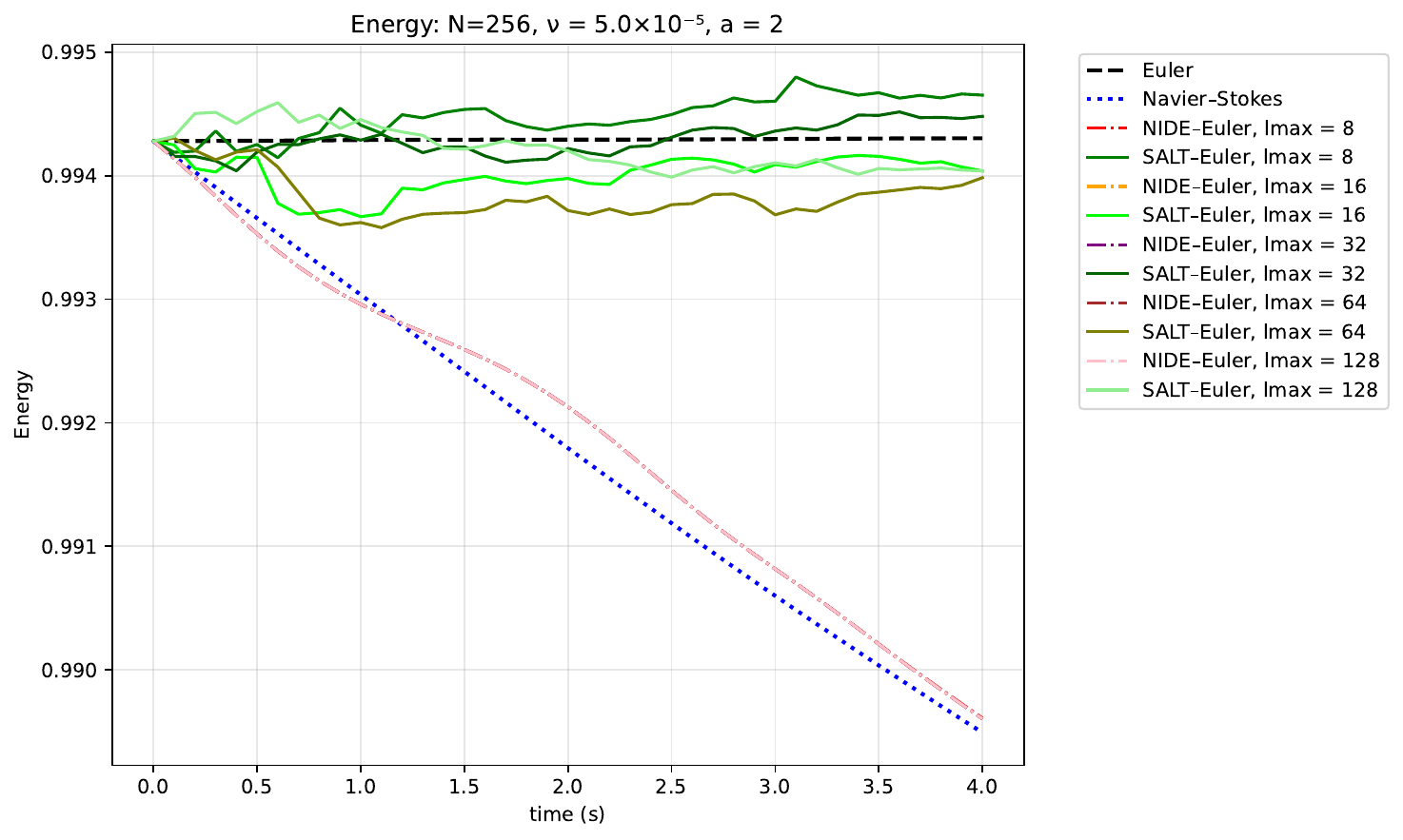}

    \includegraphics[width=0.49\linewidth]{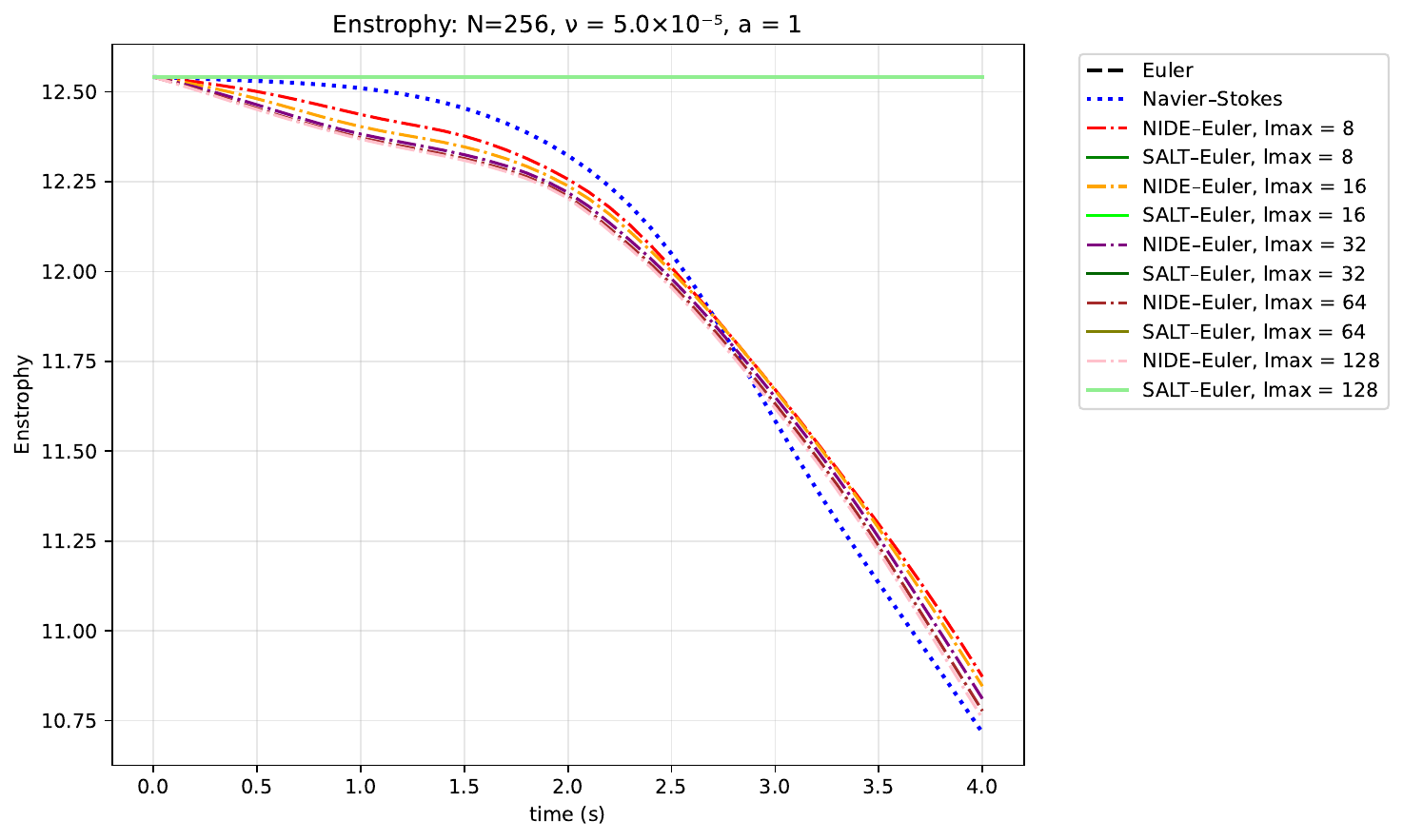}
    \includegraphics[width=0.49\linewidth]{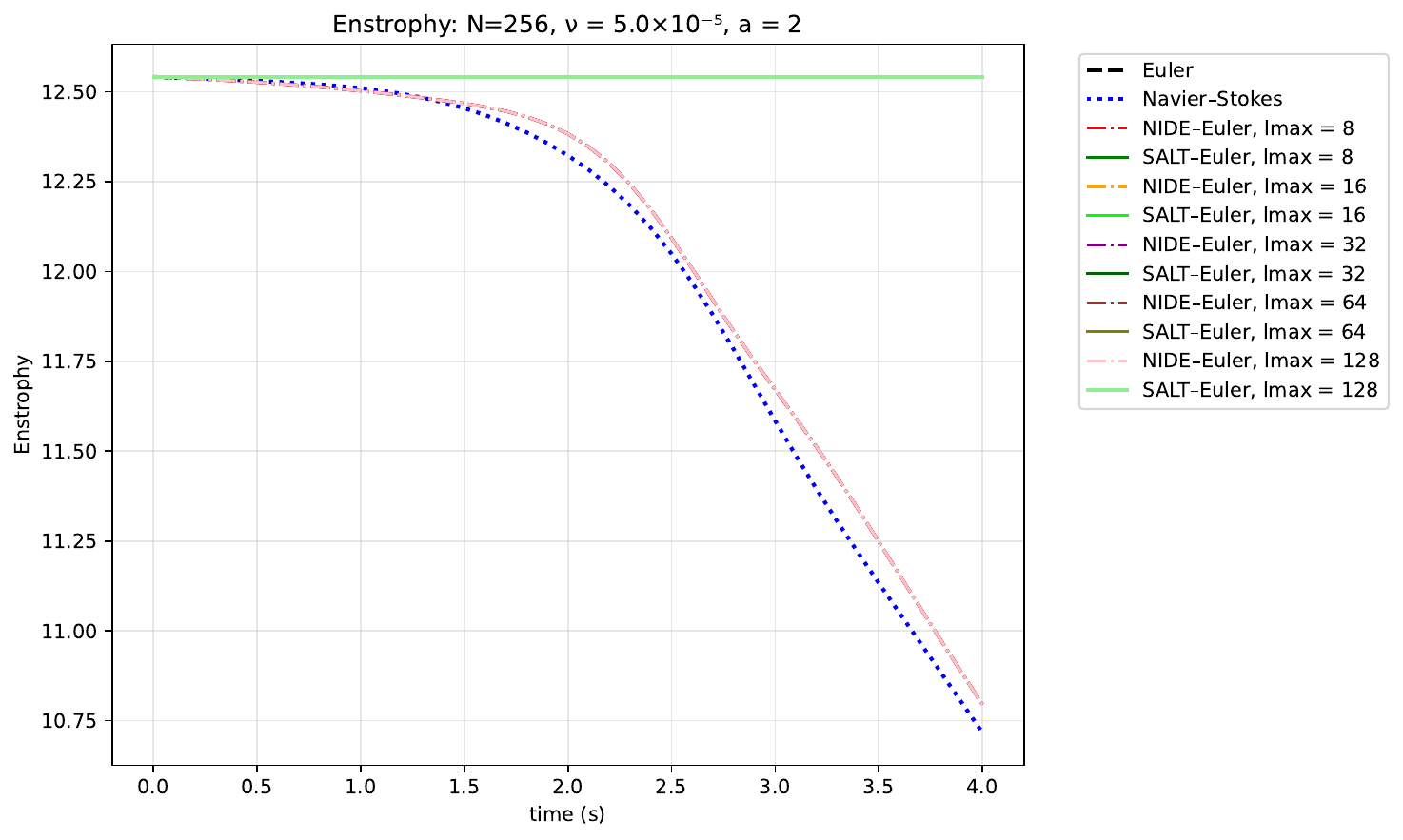}
    
    \caption{Energy evolution (top row) and enstrophy evolution (bottom row) for various deterministic and stochastic simulations at resolution $N=256$. The deterministic simulations include the standard Euler equations, the Navier--Stokes equations, and the NIDE Euler equations based on various scalings \eqref{eq:noise_scaling_numerical}.
    The ensemble mean of the energy (enstrophy) is shown for each stochastic ensemble. The noise scaling \eqref{eq:noise_scaling_numerical} is defined with $a=1$ (left column) and $a=2$ (right column).}
    \label{fig:energy_enstrophy}
\end{figure}

The evolution of the energy and the enstrophy is shown in \Cref{fig:energy_enstrophy}, where the scaling \eqref{eq:noise_scaling_numerical} defines the stochastic terms and the NIDE operator.
Each stochastic ensemble consists of 100 realizations, and the shown values are obtained by first computing the energy evolution for each realization and subsequently averaging over the entire ensemble.

As expected, the energy and enstrophy in the solution of the Euler equations remains constant while the Navier--Stokes equations and the NIDE Euler equations dissipate energy and enstrophy. 
At $a=1$, a distinct dissipative behavior is observed in the stochastic Euler equations and the NIDE Euler equations, apparent from the decrease in energy.
Namely, increasing the number of terms in the stochastic forcing or the NIDE operator increases the dissipation strength. More so, this result seems to go in the direction of showing an enhancement of diffusion respect the conjectured limit of the Navier--Stokes equation. Hence, suggesting a bridge to possible show numerically some form of "anomalous dissipation" as in \cite{agresti2024anomalous}. At the same time, the enstrophy of the stochastic Euler equations is preserved since these solutions remain on coadjoint orbits. 
This suggests the feasibility of dissipation on coadjoint orbits via transport noise.

Qualitatively different behavior is observed at $a=2$.
The energy and enstrophy in the deterministic NIDE Euler are virtually indistinguishable when varying the number of added terms. 
This is in line with the predictions made in \Cref{eq:energy_diss} and \Cref{eq:enstrophy_diss}, that as $a$ becomes larger, the higher modes become less importance. 
We thus observe that if $a = 2$, then only the first few modes are relevant and the behavior is strikingly close to the the Laplacian, alluding to the Laplacian also being a NIDE operator with only three terms.

Furthermore, a close agreement with the Navier--Stokes equations is observed. This is expected since increasing $a$ results in a larger weight on the low-degree spherical harmonics, hence the NIDE operator behaves more like the standard Laplace--Beltrami operator.
The dissipative behavior is not observed in the stochastic realizations, suggesting that the scaling limit will not take place under such regime and that the effect of the noise term balances our with the deterministic operator.
In fact, the energy is roughly maintained throughout the entire simulated time interval.

It can be noted, as a conclusion, that one of the possible reasons for observing such behaviors in the stochastic equation when $a=2$ could be traced to the fact that putting emphasis on lower modes in the NIDE operator, makes it closer to a Laplacian and, at the same time the noise term behaves like $\|\alpha^M_\cdot\|^2_\infty\sim \nu_\mathrm{SALT}$ in norm, in a way that counteracts the dissipation, and this could potentially explain the observed behavior.

\section{Conclusions}
In this paper, we investigated dissipative properties of transport noise on the sphere. 
The geometry of the two-dimensional Euler equations prompt the use of stochasticity understood in the Stratonovich sense. The corresponding \ito formulation, more suited for analysis, revealed the presence a noise-induced differential elliptic (NIDE) operator. This operator includes, e.g., the Laplace--Beltrami operator, and its energy- and enstrophy-dissipating properties were analyzed. 
The feasibility of reproducing scaling limits of transport noise were discussed.

Numerical experiments inspired by dissipative limits of transport noise were carried out using the structure-preserving Zeitlin approach for the Euler equations on the sphere. These experiments were performed using a recently developed scalable integrator for stochastic Lie--Poisson dynamics.
While it is generally infeasible to numerically compute such dissipative noise limits, the presented results indicated that dissipative behavior can be induced by transport noise with appropriately scaled parameters. Simultaneously, the stochastic realizations maintained the geometric properties also present in the solution of the deterministic Euler equation.

Future work is dedicated to further analysis of transport noise on the sphere. Specifically, improving the understanding of the NIDE operator and deriving noise scalings to achieve dissipative limits can contribute to the calibration of transport noise as a parametrization for unresolved processes in geophysical fluid simulations.

\bibliographystyle{abbrv}
\bibliography{bib}

\begin{thebibliography}{10}

\bibitem{agresti2024anomalous}
A.~Agresti.
\newblock On anomalous dissipation induced by transport noise.
\newblock {\em arXiv preprint arXiv:2405.03525}, 2024.

\bibitem{bordemann1991gl}
M.~Bordemann, J.~Hoppe, P.~Schaller, and M.~Schlichenmaier.
\newblock $\mathfrak{gl}(\infty)$ and geometric quantization.
\newblock {\em Comm. Math. Phys.}, 138:209--244, 1991.

\bibitem{bordemann1994toeplitz}
M.~Bordemann, E.~Meinrenken, and M.~Schlichenmaier.
\newblock Toeplitz quantization of {K}{\"a}hler manifolds ang gl (n), {N}→$\infty$ limits.
\newblock {\em Comm. Math. Phys.}, 165:281--296, 1994.

\bibitem{carigi2023dissipation}
G.~Carigi and E.~Luongo.
\newblock Dissipation properties of transport noise in the two-layer quasi-geostrophic model.
\newblock {\em J. Math. Fluid Mech.}, 25(2):28, 2023.

\bibitem{cotter2019numerically}
C.~Cotter, D.~Crisan, D.~D. Holm, W.~Pan, and I.~Shevchenko.
\newblock Numerically modeling stochastic {L}ie transport in fluid dynamics.
\newblock {\em Multiscale Model. Simul.}, 17(1):192--232, 2019.

\bibitem{cotter2020particle}
C.~Cotter, D.~Crisan, D.~D. Holm, W.~Pan, and I.~Shevchenko.
\newblock A particle filter for stochastic advection by lie transport: a case study for the damped and forced incompressible two-dimensional {E}uler equation.
\newblock {\em SIAM/ASA J. Uncertain. Quantif.}, 8(4):1446--1492, 2020.

\bibitem{crisan2019solution}
F.~F. Crisan~D. and H.~D.D.
\newblock Solution properties of a {3D} stochastic {E}uler fluid equation.
\newblock {\em J. Nonlinear Sci.}, 29(3):813--870, 2019.

\bibitem{dufee2022stochastic}
B.~Duf{\'e}e, E.~M{\'e}min, and D.~Crisan.
\newblock Stochastic parametrization: An alternative to inflation in ensemble {K}alman filters.
\newblock {\em Q. J. R. Meteorol. Soc.}, 148(744):1075--1091, 2022.

\bibitem{emery2006two}
M.~{\'E}mery.
\newblock On two transfer principles in stochastic differential geometry.
\newblock In {\em S{\'e}minaire de Probabilit{\'e}s XXIV 1988/89}, pages 407--441. Springer, 2006.

\bibitem{ephrati2024exponential}
S.~Ephrati, E.~Jansson, A.~Lang, and E.~Luesink.
\newblock An exponential map free implicit midpoint method for stochastic {L}ie-{P}oisson systems.
\newblock {\em arXiv preprint arXiv:2408.16701}, 2024.

\bibitem{ephrati2023data}
S.~R. Ephrati, P.~Cifani, E.~Luesink, and B.~J. Geurts.
\newblock Data-driven stochastic {L}ie transport modeling of the 2{D} {E}uler equations.
\newblock {\em J. Adv. Model. Earth Syst.}, 15(1):e2022MS003268, 2023.

\bibitem{flandoli2021scaling}
F.~Flandoli, L.~Galeati, and D.~Luo.
\newblock Scaling limit of stochastic {2D} {E}uler equations with transport noises to the deterministic {N}avier--{S}tokes equations.
\newblock {\em J. Evol. Equ.}, 21(1):567--600, 2021.

\bibitem{FlaLuo1}
F.~Flandoli and D.~Luo.
\newblock {Convergence of transport noise to {O}rnstein–{U}hlenbeck for 2{D} {E}uler equations under the enstrophy measure}.
\newblock {\em Ann. Probab.}, 48(1):264 -- 295, 2020.

\bibitem{FlaPapMor}
F.~Flandoli, S.~Morlacchi, and A.~Papini.
\newblock Effect of transport noise on {K}elvin--{H}elmholtz instability.
\newblock In B.~Chapron, D.~Crisan, D.~Holm, E.~M{\'e}min, and A.~Radomska, editors, {\em Stochastic Transport in Upper Ocean Dynamics II}, pages 29--52, Cham, 2024. Springer Nature Switzerland.

\bibitem{flandoli2023average}
F.~Flandoli, A.~Papini, and M.~Rehmeier.
\newblock Average dissipation for stochastic transport equations with l{\'e}vy noise.
\newblock In {\em Stochastic Transport in Upper Ocean Dynamics Annual Workshop}, pages 45--59. Springer Nature Switzerland Cham, 2023.

\bibitem{GayBalmaz2013}
F.~Gay-Balmaz and D.~D. Holm.
\newblock Selective decay by casimir dissipation in inviscid fluids.
\newblock {\em Nonlinearity}, 26(2):495–524, Jan. 2013.

\bibitem{goodair20223d}
D.~Goodair and D.~Crisan.
\newblock On the {3D} navier-stokes equations with stochastic lie transport.
\newblock In {\em Stochastic Transport in Upper Ocean Dynamics Annual Workshop}, pages 53--110. Springer Nature Switzerland Cham, 2022.

\bibitem{holm2015variational}
D.~D. Holm.
\newblock Variational principles for stochastic fluid dynamics.
\newblock {\em Proc. R. Soc. Lond. Ser. A Math. Phys. Eng. Sci.}, 471(2176):20140963, 2015.

\bibitem{hoppe1989diffeomorphism}
J.~Hoppe.
\newblock Diffeomorphism groups, quantization, and {SU}($\infty$).
\newblock {\em Internat. J. Mod. Phys. A}, 4(19):5235--5248, 1989.

\bibitem{hoppe1998some}
J.~Hoppe and S.-T. Yau.
\newblock Some properties of matrix harmonics on $\mathbb{S}^2$.
\newblock {\em Comm. Math. Phys.}, 195:67--77, 1998.

\bibitem{kloeden1992stochastic}
P.~E. Kloeden, E.~Platen, P.~E. Kloeden, and E.~Platen.
\newblock {\em Stochastic differential equations}.
\newblock Springer, 1992.

\bibitem{lang2023well}
O.~Lang and D.~Crisan.
\newblock Well-posedness for a stochastic {2D} euler equation with transport noise.
\newblock {\em Stoch. Partial Differ.}, 11(2):433--480, 2023.

\bibitem{memin2014fluid}
E.~M{\'e}min.
\newblock Fluid flow dynamics under location uncertainty.
\newblock {\em Geophys. Astrophys. Fluid Dyn.}, 108(2):119--146, 2014.

\bibitem{milstein2004stochastic}
G.~N. Milstein and M.~V. Tretyakov.
\newblock {\em Stochastic numerics for mathematical physics}, volume~39.
\newblock Springer, 2004.

\bibitem{modin2020casimir}
K.~Modin and M.~Viviani.
\newblock A {C}asimir preserving scheme for long-time simulation of spherical ideal hydrodynamics.
\newblock {\em J. Fluid Mech.}, 884:A22, 2020.

\bibitem{modin2020lie}
K.~Modin and M.~Viviani.
\newblock Lie--{P}oisson methods for isospectral flows.
\newblock {\em Found. Comput. Math.}, 20(4):889--921, 2020.

\bibitem{modin2024two}
K.~Modin and M.~Viviani.
\newblock Two-dimensional fluids via matrix hydrodynamics.
\newblock {\em arXiv preprint arXiv:2405.14282}, 2024.

\bibitem{resseguier2017geophysical}
V.~Resseguier, E.~M{\'e}min, and B.~Chapron.
\newblock Geophysical flows under location uncertainty, part {II} {Q}uasi-geostrophy and efficient ensemble spreading.
\newblock {\em Geophys. Astrophys. Fluid Dyn.}, 111(3):177--208, 2017.

\bibitem{Sadourny1985}
R.~Sadourny and C.~Basdevant.
\newblock Parameterization of subgrid scale barotropic and baroclinic eddies in quasi-geostrophic models: Anticipated potential vorticity method.
\newblock {\em J. Atmos. Sci.}, 42(13):1353–1363, July 1985.

\bibitem{zeitlin2004self}
V.~Zeitlin.
\newblock Self-consistent finite-mode approximations for the hydrodynamics of an incompressible fluid on nonrotating and rotating spheres.
\newblock {\em Phys. Rev. Lett.}, 93(26):264501, 2004.

\end{thebibliography}

\end{document}